\documentclass[twocolumn,showpacs,preprintnumbers,amsmath,amssymb]{revtex4-1}
\usepackage[dvips]{graphicx}
\usepackage{epstopdf}
\usepackage{float}
\DeclareMathOperator*{\fmin}{min}

\newcommand{\bec}{\begin{center}}
\newcommand{\enc}{\end{center}}
\newcommand{\be}{\begin{equation}}
\newcommand{\ee}{\end{equation}}
\newcommand{\bmi}{\begin{minipage}}
\newcommand{\emi}{\end{minipage}}

\newcommand{\bi}{\begin{itemize}}
\newcommand{\ei}{\end{itemize}}
\newcommand{\ba}{\begin{array}}
\newcommand{\ea}{\end{array}}

\newcommand{\lqu}{\left[}
\newcommand{\rqu}{\right]}

\newcommand{\ora}[1]{\overrightarrow{#1}}

\begin{document}
\title{A discrete-pulse optimal control algorithm with an application to
spin systems}
\author{G. Dridi$^{1,2,3}$}
\author{M. Lapert$^4$}
\author{J. Salomon$^5$}
\author{S. J. Glaser$^4$}
\author{D. Sugny$^{3,6}$}
\email{dominique.sugny@u-bourgogne.fr}

\affiliation{$^1$ Nanosciences Group and MANA Satellite, CEMES-CNRS, 29 Rue Jeanne Marvig, F- 31055, Toulouse, France}
\affiliation{$^2$ CMAP, UMR 7641, Ecole Polytechnique CNRS, Route de Saclay 91128 Palaiseau Cedex France}
\affiliation{$^3$ Laboratoire Interdisciplinaire Carnot de
Bourgogne (ICB), UMR 6303 CNRS-Universit\'e de Bourgogne, 9 Av. A.
Savary, BP 47 870, F-21078 DIJON Cedex, FRANCE}
\affiliation{$^4$ Department of Chemistry, Technische Universit\"at
M\"unchen, Lichtenbergstrasse 4, D-85747 Garching, Germany}
\affiliation{$^5$ CEREMADE,
Universit\'e Paris Dauphine, Place du Mar\'echal De Lattre De Tassigny, 75775 Paris Cedex 16, France}
\affiliation{$^6$ Institute for Advanced Study, Technische Universit\"at M\"unchen, Lichtenbergstrasse 2 a, D-85748 Garching, Germany}
\date{\today}
\begin{abstract}
This article is aimed at extending the framework of optimal
control techniques to the situation where the control field values are restricted to a finite set. We propose a
generalization of the standard GRAPE algorithm suited to this
constraint. We test the validity and the efficiency of this approach
for the inversion of an inhomogeneous ensemble of spin systems with different offset frequencies. It
is shown that a remarkable efficiency can be achieved even for a
very limited number of discrete values. Some applications in Nuclear
Magnetic Resonance are discussed.
\end{abstract}
 \maketitle
\section{Introduction}
\begin{figure}[htpb]
\includegraphics[scale=1.25]{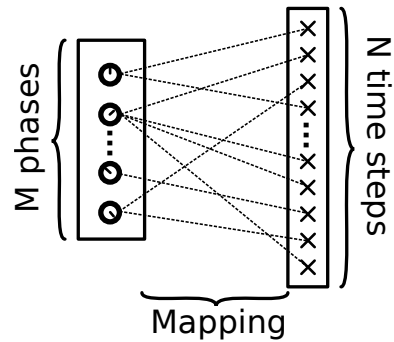}
\caption{Schematic representation of the mapping vector between the set of $M$- phases and the set of $N$- time steps. The mapping is depicted by the dashed lines.}\label{fig:sketch}
\end{figure}

The design of control sequences accounting for experimental constraints is a central task in a variety of domains in quantum dynamics extending from photochemistry, Nuclear Magnetic Resonance (NMR) and quantum information science \cite{cat,reviewQC,warren,rice,shapiro,tannorbook,spin,garon,imaging}. Nowadays, Optimal Control Theory (OCT) reveals to be a highly efficient and versatile tool to bring answers to the different issues raised by the experimental setups \cite{reviewQC,pont,krotov,bryson,koch,lapertglaser,simutime,khaneja,kosloff,somloi,zhu,grapeino,maday,salomon2}. For the past few years, there has been an intense theoretical activity in developing new optimal control procedures able to build high quality control fields in presence of some experimental imperfections and constraints \cite{reviewQC,koch,gross,glaseramp,grapetheory}. These include spectral constraints \cite{viviespectrum,lapertspectrum,palao,glaserspectrum}, amplitude and  phase transients \cite{Prisner_2012}, non-linear interactions between the system and the control field \cite{ohtsukinonlinear,lapert1,zhang}, robustness against experimental uncertainties and errors \cite{grapeino,daems}. To date, the majority of studies has assumed that the amplitude and phase of the control field can vary continuously. For example, in modern NMR spectrometers \cite{spin}, the amplitude and phase of the control pulses can be defined with high resolution, allowing for a virtually continuous variation of these parameters  \cite{Keeler_1992,Sternin_1995}.
More generally, this is possible in experimental settings, where arbitrary waveform generators are available \cite{Prisner_2012,Martinis_2013,Wunderlich_2013,Bowler_2013}.
However, in many cases the available hardware only allows
to switch between a discrete set of pulse phases \cite{Haeberlen_1996, Smith_2009}. This set of phases can be chosen before the experiment, but cannot be altered during the experiment.
Hence the control is quantized and restricted to a fixed finite number of values, which can nevertheless be optimized (See Fig.~\ref{fig:sketch} for a schematic description).

This paper is aimed at tackling this problem by proposing two fundamentally different approaches. The first method is based on an extension of a standard optimal iterative procedure, namely GRAPE \cite{grapetheory,grapereview}. The second option can be viewed as a brute force strategy in the sense that a standard optimization method is used in a first step, producing a continuous control field. A Lloyd procedure, analogous in its spirit to a mean square method, is then applied to quantize the field \cite{lloyd}. We test the validity and the efficiency of the two approaches on a benchmark control problem, the simultaneous control of an ensemble of uncoupled spin 1/2 particles with different offset terms \cite{spin,inversion,inversion2}. Extensive numerical simulations reveal that an efficient control can be achieved even for a small number of quantized values of the control field. We analyze the relative efficiency of the two algorithms as a function of the number of allowed values for the control field. Finally, note that the discretization of the control field has been recently used in a series of paper to accelerate the numerical integration of the time-dependent Schr\"odinger equation \cite{toolkitnew1,toolkitnew2,toolkitnew3,toolkitnew4,toolkitnew5}. This idea based on the precomputation of elementary propagators was applied to quantum optimal control problem and would be also useful in our case.

The remainder of the paper is organized as follows. In Sec. \ref{sec2}, we outline the principles of the new optimization procedures, paying special attention to the flexibility and to the applicability of the methods. Section \ref{sec3} is dedicated to the presentation of the numerical results. The efficiency of the quantized control field is compared to its continuous counterpart. The advantages of the optimization algorithm directly accounting for the constraint and of the brute force strategy are also discussed. Conclusion and prospective views are given in Sec. \ref{sec4}.
\section{Theory}\label{sec2}
The goal of this section is to propose an optimal control
algorithm suited to quantum systems where the control field is restricted to a finite set $\mathcal{M}$ of $M$ discrete values
$\mathcal{M}=\{v_1,v_2,\cdots,v_M\}$, where $M$ is fixed and the
values $v_k$ are optimized. To simplify the discussion, we
consider here the case of mixed quantum systems, which also applies directly to pure states and can be straightforwardly extended to the control of unitary operators in quantum computing. Starting from an initial state $\rho_0$, optimal control is invoked in order to
maximize the projection onto a target state $\rho_f$, which is
measured by the following figure of merit:
\begin{equation}\label{eq1}
\Phi=\textrm{Tr}[\rho_f^\dag\rho(t_f)],
\end{equation}
where $t_f$ is the control duration and $\rho(t)$ the state of the
system at time $t$.
\subsection{The GRAPE algorithm}\label{sec2a}
In this paper, we consider specifically the GRAPE algorithm
\cite{grapetheory}, but the same construction of the discrete
version can be used for other algorithms such as the monotonic or Krotov
ones \cite{tannor,ohtsuki,salomon,bifurcating1,lapert1} (see also the general analysis of such methods \cite{toolkit1,toolkit2}).

The time evolution of the mixed-state $\rho(t)$ is ruled by the
Liouville-von Neumann equation:
\begin{equation}\label{eq2}
i\dot{\rho}(t)=[H(t),\rho(t)],
\end{equation}
where units such that $\hbar=1$ have been chosen. The Hamiltonian
$H(t)$ of the system can be written as:
\begin{equation}\label{eq3}
H(t)=H_0+u(t)H_1,
\end{equation}
where $H_0$ is the field-free Hamiltonian operator and $H_1$ the
interaction part. 
A general presentation of the algorithm is proposed here, but the procedure can be straightforwardly adapted to phase-modulated pulses with a constant amplitude \cite{phaseglaser}, as used in Sec.~\ref{sec3}.

Let
$U(t)$ be the evolution operator associated with the Hamiltonian
$H(t)$, such that $\rho(t)=U(t)\rho_0U^\dag (t)$. Following the description of a standard GRAPE algorithm \cite{grapetheory}, we assume that the field
$u$ is discretized in time with a time step $\Delta t=t_f/N$, where
$N$ is the number of values of the field. We denote by $u_j$, with $j=1,2,\cdots, N$ the value of $u$ in the
interval $[(j-1)\Delta t,j\Delta t]$. In this approximation, note
that the control field $u$ is now described by a set of $N$ reals
$(u_1,u_2,\cdots, u_N)$. The time evolution can be computed by the $N$ block propagators $(U_1,U_2,\cdots,U_N)$,
where
\begin{equation}\label{eq4}
U_j=e^{-i\Delta t (H_0+u_jH_1)}.
\end{equation}
The optimal control problem is described by a figure of merit $\Phi$ to be maximized. Using the different evolution operators $U_j$, $\Phi$ can be written as follows:
\begin{equation}\label{eq7}
\Phi=\textrm{Tr}\big(U_NU_{N-1}\cdots U_1\rho_0 U_1^\dag U_2^\dag
\cdots U_N^\dag  \rho_f \big).
\end{equation}
The optimization procedure is based on the derivation of the gradient of the figure of merit with respect to the different values of the control field:
\begin{equation}
\begin{array}{rl}
\frac{\partial \Phi}{\partial u_j}=&\textrm{Tr}\big(U_NU_{N-1}\cdots \frac{\partial U_j}{\partial u_j}\cdots U_1\rho_0 U_1^\dag U_2^\dag \cdots  U_N^\dag \rho_f \big)\\
+&\textrm{Tr}\big(U_NU_{N-1}\cdots U_1\rho_0 U_1^\dag U_2^\dag \cdots \frac{\partial U_j^\dag}{\partial u_j}\cdots U_N^\dag \rho_f \big).
\end{array}
\end{equation}
We assume that the different propagators
are approximated by a first order split-operator, which is valid
up to the order 2 in $\Delta t$:
\begin{equation}\label{eq5}
U_j\simeq e^{-i\Delta t H_0}e^{-i\Delta t u_jH_1}.
\end{equation}
This hypothesis simplifies the computation of the derivative of $U_j$
with respect to $u_j$, which can be written as:
\begin{equation}\label{eq6}
\frac{\partial U_j}{\partial u_j}\simeq U_j(-i\Delta t H_1).
\end{equation}
Note that a different approximation of the propagator would give a different derivative. A straightforward computation using Eq.~(\ref{eq6}) then leads to the gradient of the figure of merit:
\begin{equation}\label{eq8}
\frac{\partial \Phi}{\partial u_j}\simeq -i\Delta t
\textrm{Tr}\big(\lambda_j [H_1,\rho_j]\big),
\end{equation}
where the states $\rho_j$ and $\lambda_j$ are defined by:
\begin{eqnarray*}
& & \rho_j=U_{j-1}\cdots U_1 \rho_0 U_1^\dag \cdots U_{j-1}^\dag, \\
& & \lambda_j=U_j^\dag U_{j+1}^\dag \cdots U_N^\dag \rho_f U_N
\cdots U_j.
\end{eqnarray*}
In a first-order GRAPE algorithm (see \cite{grapesecond} for a
recent generalization to the second order), the field $u$ is updated at each step by the
formula:
\begin{equation}\label{eq9}
u_j\mapsto u_j+\varepsilon \frac{\partial \Phi}{\partial u_j},
\end{equation}
with $\varepsilon$ a small real parameter, which is chosen from a line search method to ensure the increase of the figure of merit $\Phi$.
\subsection{A discrete-pulse version of GRAPE}\label{sec2b}
In the continuous version of the GRAPE algorithm, the control
field $u$ can take any real value, \textit{i.e.} $u_j\in\mathbb{R}$. We
consider in this section that $u$ is restricted to a finite set
$\mathcal{M}$ of $M$ values. We introduce the mapping vector $\ora{p}$
from the set $\{1,2,\cdots,N\}$ to $\{1,2,\cdots M\}$ which
associates with the different values of $u$, a value in
$\mathcal{M}$: $u_j=v_{\ora{p}(j)}$. This mapping makes a correspondence between the set of time steps and the set of dicrete values. A schematic illustration of this process is given in Fig.~\ref{fig:sketch}.

This optimization procedure has a non trivial character in the sense that both the discrete values $\{v_m\}$ and the mapping $\ora{p}$ can be adjusted to maximize the figure of merit $\Phi$. In the algorithm proposed in this paper, each iteration is decomposed into two sub-steps.
In the first sub-step, the mapping $\ora{p}$ is fixed and the $\{v_m\}$- values are optimized through a modified version of GRAPE which can be described as follows. The functional $\Phi$ can be derived with respect to
$v_m$:
\begin{equation}\label{eq11}
\frac{\partial \Phi}{\partial v_m}\simeq \sum_{j | \ora{p}(j)=m}-i\Delta t
\textrm{Tr}\big(\lambda_j [H_1,\rho_j]\big).
\end{equation}
Note that the derivative of the quality factor with respect to the discrete value $v_m$ is simply the sum over the derivatives with respect to $u_j$ which maps to $v_m$.
The control $v_m$ is then updated at each step of the algorithm through the
formula:
\begin{equation}\label{eq12}
v_m\mapsto v_m+\varepsilon \frac{\partial \Phi}{\partial v_m}.
\end{equation}
 The roles are reversed in the second stage which aims at optimizing $\ora{p}$, while the discrete values $\{v_m\}$ are not changed. Since the total number $N$ of values of the field can be very large, it would be time consuming to find the global optimal mapping $\ora{p}$. Instead, we adopt a more direct approach which allows to improve the final result without no guarantee to attain its upper bound. The efficiency of this procedure will be shown numerically and discussed in Sec.~\ref{sec3}.

For each time step $j\in \{1,2,\cdots,N\}$ taken in increasing order, we test the $M$ possible values of the control field $u_j=v_l$, $l=1,2,\cdots,M$, by computing the corresponding figure of merit $\Phi(u_j=v_l)$. Then we define the new discrete phase as being the one which maximizes the quality factor:
$$
  \ora{p}(j) = k = \arg\max \Phi(v_k)
$$
A proper use of the adjoint state allows us to obtain a computational cost of the algorithm with depends linearly on $M$.

Since our method is a two-step procedure, it cannot be easily extended to the second order (the second order is related here to the gradient and not to the order of accuracy of the propagator). Note also that this approach is compatible with a toolkit method \cite{toolkitnew1,toolkitnew2,toolkitnew3,toolkitnew4,toolkitnew5}. At each step of the algorithm, the propagators associated with the set of $M$- discrete values can be precomputed. These propagators can be used to update the mapping at iteration $k-1$ and the same set of propagators allows us to derive the gradient at step $k$ of the algorithm.
\subsection{Quantization and Lloyd's algorithm}\label{sec2c}
We propose in this paragraph a second strategy based on Lloyd's algorithm. The idea consists first in using
a standard GRAPE algorithm to build a continuous control field. In our numerical example, we will be interested in a case where the phase of the field is optimized, while its amplitude is constant. Nevertheless, it would be straightforward to extend this procedure to a general situation with no constraint on the control field. Lloyd's approach is then applied in a second step to quantize this field. Lloyd's algorithm, which is explained in detail below (we refer the reader to \cite{lloyd} for additional information), can be viewed as a mean-square approximation procedure which allows us to approach a continuous function by a discrete one. In contrast to the discrete version of GRAPE presented in Sec. \ref{sec2b}, no information about the dynamics is used for the quantization. In other words, the computation of the continuous control field is sufficient to design its discrete counterpart. However, since only a geometric (and not a dynamical) information is used, there is a priori no guarantee about the efficiency of the quantized control.

To simplify the description of the algorithm, we consider a set $(u_i)_{i=1,\ldots,N}$ of values in the interval $[0,2\pi[$ (see Sec. \ref{sec3} for details), which correspond to the pulse phases derived from the continuous GRAPE algorithm. The angles $(u_i)$ are defined modulo $2\pi$ and we identify 0 and $2\pi$ by periodicity. The discrete control field values derived from Lloyd's algorithm will be denoted $(\omega_i)_{i=1,\ldots,N}$. We also introduce two sets of $M$ reals belonging to $[0,2\pi[$, $B^{(k)}=(B_j^{(k)})_{j=1,\ldots,M}$ and $Y^{(k)}=(Y_j^{(k)})_{j=1,\ldots,M}$ sorted in increasing order, $k$ being the iteration step of the algorithm. The initial set $B^{(0)}$ is defined such that the values $B_j^{(0)}$ are equally distributed in increasing order in the interval $[0,2\pi[$. Note that the $2\pi$- periodicity of the control field is taken into account in the different relations used in the algorithm, even if this point is not explicitly mentioned below in order to clarify the presentation of the computation of the discrete field. The algorithm can be described as follows.



\begin{figure}[htpb]
\includegraphics[scale=0.75]{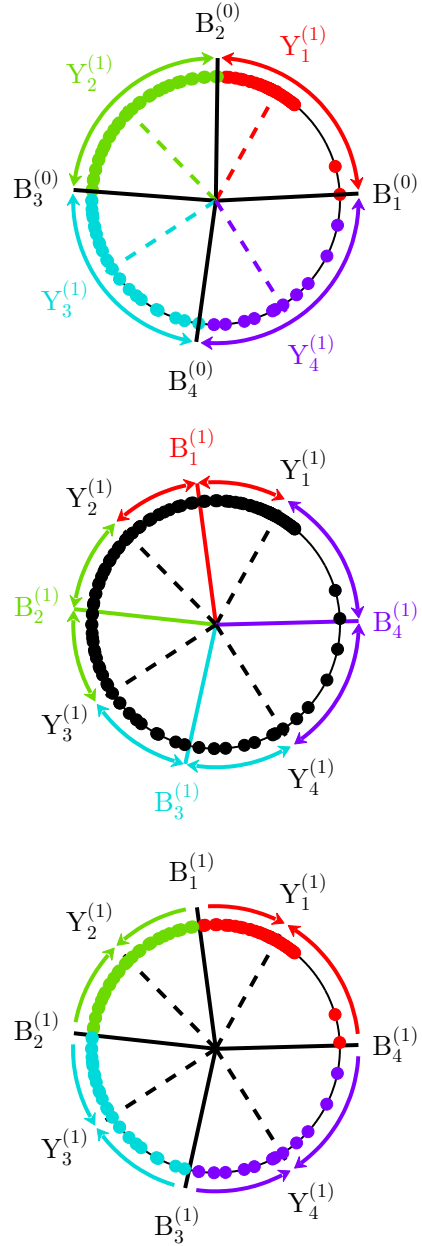}
\caption{(Color online) Schematic representation of the different steps of Lloyd's algorithm in the case $M=4$. The dots indicate the position of the control field values $\{u_i\}$. The two sets of reals $B^{(k)}$ and $Y^{(k)}$ are represented by solid and dashed lines, respectively. At iteration $k=1$, the top panel shows the way to compute the $Y^{(1)}_j$ values according to the step 2 of the algorithm. Different colors have been used to help the understanding of the procedure. Starting from this new set $Y^{(1)}$, the set $B^{(1)}$ is defined in the middle panel by using step 4 of the algorithm. The bottom panel depicts a new step 2 for the next iteration of the algorithm. At this stage, note that the new $Y$- values have not been computed.}\label{fig1}
\end{figure}

{\bf Algorithm:} Given the initial set $B^{(0)}$, the control phases $(u_i)_{i=1,\ldots,N}$, and the iteration parameters $k=0$,
$\epsilon>0$, $J_0=0$, do:
\begin{enumerate}
  \item $k=k+1$.
  \item For $j=1,\cdots,M$, do\\
  $Y_j^{(k)}={\rm Mean}_i(u_i\in\lqu B^{(k-1)}_j,B^{(k-1)}_{j+1}\rqu)$
  \item For $i=1,\cdots,N$, evaluate\\
  $d_i=\fmin_{j} |u_i-Y_j^{(k)}|~{\rm mod}~2\pi$\\
  $J_k=\sum_{i=1}^Nd_i$
  \item if $|J_k-J_{k-1}|>\epsilon$, do \\
	  For $j=1,\cdots,M$, $B_j^{(k)}=(Y_j^{(k)}+Y_{j+1}^{(k)})/2$\\
	    and go to 1\\
	else go to 5\\
  \item For $i=1,\cdots,N$, compute $\omega_i={\rm Proj}(u_i,Y^{(k)})$
\end{enumerate}
Note that the function  "${\rm Mean}$" stands for the mean value of a set of numbers belonging to a given interval. The function "${\rm Proj}$" denotes the projection of each $u_i$ onto the discrete set of values $Y^{(k)}$. More precisely, we define the projection as the value $\omega_i\in Y^{(k)}$ which minimizes the distance from $u_i$ to $Y^{(k)}$. It is then straightforward to define the mapping $\ora{p}$ and the set $\mathcal{M}$ from the values $(\omega_i)_{i=1,\ldots,N}$ such that $v_{\ora{p}(j)}=\omega_i$. A schematic description of this algorithm is displayed in Fig.~\ref{fig1}.

\section{Application to the control of spin systems}\label{sec3}
This section is dedicated to the application of the
discrete version of the GRAPE algorithm for controlling the dynamics of
spin systems driven by radio-frequency magnetic fields. The different numerical values are chosen so as to reproduce the typical features of a spin sample in
liquid state NMR spectroscopy with $B_0$- inhomogeneities \cite{spin}. We test the validity and the efficiency of the two approaches developed above for the inversion of an ensemble of inhomogeneous uncoupled spins with different resonant offset frequencies \cite{inversion}.
\subsection{The model system}

In a given rotating frame, the Hamiltonian of each isochromat, \textit{i.e.} each sub-system with a different resonance frequency $\omega$, is $H_{\omega}=\omega\sigma_z + \omega_x\sigma_x + \omega_y\sigma_y$. For a matter of convenience and to simplify the numerical implementation of this problem, we move to the Bloch picture \cite{spin}:
\begin{align}\label{eq13}
\dot{M}_{x}^{(\omega)}=&-\omega M_{y}+\omega_{y}M_{z}\\ \nonumber
\dot{M}_{y}^{(\omega)}=& \omega M_{x}-\omega_{x}M_{x}\\
\dot{M}_{z}^{(\omega)}=&\omega_{x} M_{y}-\omega_{y}M_{x},\nonumber
\end{align}
where the Bloch vector of the isochromat $\vec{M}^{(\omega)}=(M_x^{(\omega)},M_y^{(\omega)},M_z^{(\omega)})$ represents the state of the system, which can be readily related to the corresponding density matrix through the identification $M_i^{\omega}=\textrm{Tr}(\rho^{\omega}\sigma_i)$, with $i=\{x,y,z\}$ and $\sigma_i$ the Pauli matrices. On the right hand side of Eq. (\ref{eq13}), $\omega$ is the offset frequency term which belongs to the interval $\left[-\omega_{max}, \omega_{max}\right]$. The two components of the control field along the $x$- and $y$- directions are denoted $\omega_{x}$ and $\omega_{y}$, respectively. Starting from $\overrightarrow{M}(0)$=$\overrightarrow{M}_{z}$, the goal of the control is to reach the south pole of the Bloch sphere $\overrightarrow{F}=-\overrightarrow{M}_{z}$ for any spin of the ensemble. The quality factor or the figure of merit to maximize can be written as follows:
\begin{equation}
\Phi=\frac{1}{n_{off}}\sum^{n_{off}}_{i=1}\overrightarrow{M}^{(\omega_i)}(t_{f}).\overrightarrow{F},
\end{equation}
where $n_{off}$ is the total number of uncoupled spins and $t_{f}$ the total control time. The offset terms $\omega_i$ are chosen equally spaced in $\left[-\omega_{max}, \omega_{max}\right]$. Note that this control problem has been extensively investigated in the literature in the standard situation where the field is continuous. We refer the interested reader to key publications on this subject \cite{grapeino}. Following Ref. \cite{inversion}, the numerical values are taken to be $\omega_{max}/(2\pi)=10~\text{kHz}$, $n_{off}=200$ and $t_{f}=0.18~\text{ms}$. The control field $\overrightarrow{\omega}=\left(\omega_{x}, \omega_{y}\right)$ is assumed of fixed control amplitude $\omega_0$ and can be expressed as
\begin{equation}
\overrightarrow{\omega}=\omega_{0}\left[\cos(\theta(t))\overrightarrow{x}+\cos(\theta(t))\overrightarrow{y}\right],
\end{equation}
where $\theta\in [0,2\pi[$ is the phase to optimize. The maximum pulse amplitude $\omega_{0}/(2\pi)$ is chosen constant and equal to 10 kHz, the same value as $\omega_{max}$. The time digitization is taken as $\Delta t$=0.5 $\mu$s. We have checked that the qualitative conclusions of this paper do not depend on a specific choice of the used constants.
\subsection{Numerical results}
We are now in a position to check the efficiency of the
discrete GRAPE algorithm through the comparison with the
continuous optimal solution. In all the numerical simulations with the continuous version, we choose for the initial control phase the following simple form:
\begin{equation}
\theta(t)=\frac{\pi}{2}\left(\frac{2t}{t_{f}}-1\right)^{2}.
\end{equation}
This parabolic behavior leads to a linear evolution of the corresponding frequency. This class of control fields is known to be robust to experimental imperfections in the adiabatic limit \cite{vitanov}. With this initial guess, the GRAPE algorithm converges to the target state with an accuracy better than 0.9982 and with an optimal control field very close to the solution proposed in Ref. \cite{inversion}, see also Fig. \ref{fig5} for a plot of this optimal field.

\begin{figure}[htpb]
\includegraphics[scale=0.6]{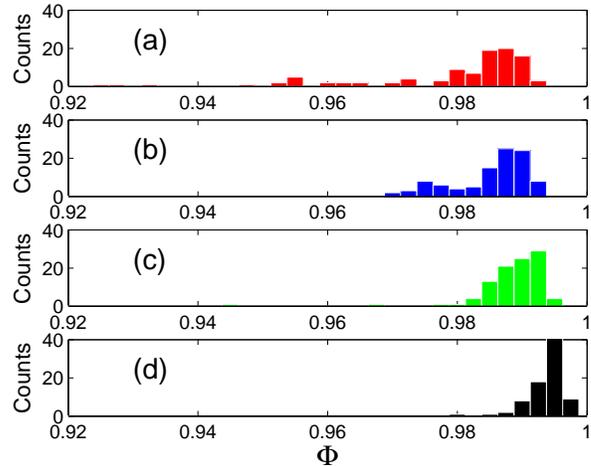}
\caption{Histogram distribution of the final quality factor $\Phi$ obtained with the discrete GRAPE method, for the cases $M=4$ ((a)- red bars), $M=8$ ((b)- blue bars), $M=12$ ((c)- green bars) and $M=16$ ((d)- black bars). In each optimal computation, the initial discrete values are randomly chosen in the interval $[0,2\pi[$ to generate 100 possible realizations.}\label{fig2}
\end{figure}

\begin{figure}[htpb]
\includegraphics[scale=0.6]{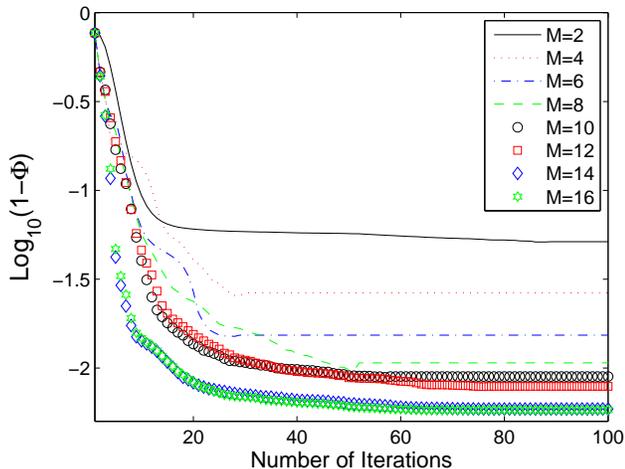}
\caption{(Color online) Plot of the quality factor $\Phi$ obtained with the discrete GRAPE method as a function of the number of iterations. For this example, the initial values of the phase are uniformly distributed in the interval $[0,2\pi[$.}\label{fig3}
\end{figure}

In the discrete case, we recall that $M$ is the number of possible discrete values that can be taken by the control phase $\theta$. Figure~\ref{fig2} shows the histogram distribution of  the quality factor $\Phi$ for the cases $M=4$, 8, 12 and 16. In each situation, the initial discrete values are randomly chosen in the interval $[0, 2\pi[$ to generate 100 possible realizations. The initial mapping $\vec{p}$ is also generated randomly. We observe that the final quality factor depends on the initially chosen (random) discrete set of values of the control field, but a significant number of examples converges towards a quality factor close to 1, even for $M=4$. For instance, Fig.~\ref{fig2} shows that a large percentage of the optimizations is close to the maximum quality factors of 0.992, 0.993, 0.994 and 0.9965 for the cases of $M=4$, $M=8$, $M=12$ and $M=16$, respectively. As could be expected, we observe that the higher the number of discrete values $M$, the narrower the distribution is. Note also that the global shape of the histogram distribution in this discrete setting is very similar to standard distributions that can be observed in the continuous case \cite{inversion,inversion2}.
Another standard choice for the initial values of the control field is an equidistant repartition in the interval $[0,2\pi[$. Here, we consider an initial zero control field and we determine the initial mapping by a forward propagation where the $M$ possible values of the field are tested for each time step. Note that this choice is not crucial for the final efficiency of the algorithm, but it allows in one shot to achieve a very high performance. The corresponding figures of merit are depicted in Fig.~\ref{fig3} as a function of the number of iterations. Here again, we observe the good behavior of the algorithm since a performance higher than 0.99 is achieved for $M\geq 8$.

\begin{figure}[htpb]
\includegraphics[scale=0.6]{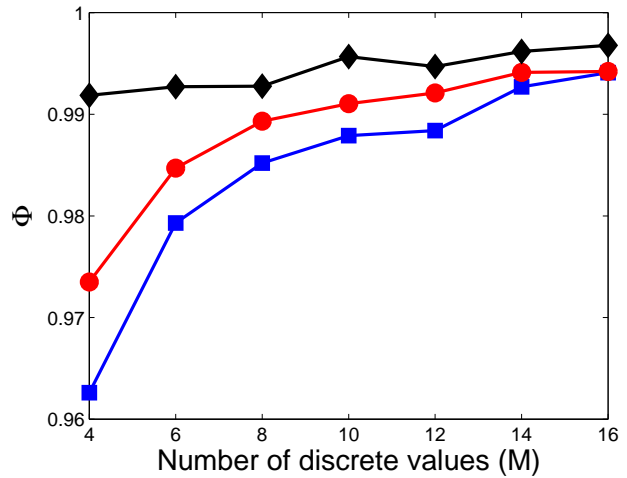}
\caption{(Color online) Plot of the final quality factor $\Phi$ for two specific initial sets of discrete values as a function of the number of discrete values $M$. In the case of a random initial set (shown by the square blue solid line), the quality factor has been obtained by averaging over the 100 realizations given in Fig.~\ref{fig2}. The red solid line with circles represents the value of the final figure of merit $\Phi$ when the initial discrete values of the control field are uniformly distributed in the interval $[0,2\pi[$. The black curve with diamonds displays the highest quality factor among the 100 random realizations.}\label{fig4}
\end{figure}

\begin{figure}[htpb]
\includegraphics[scale=0.65]{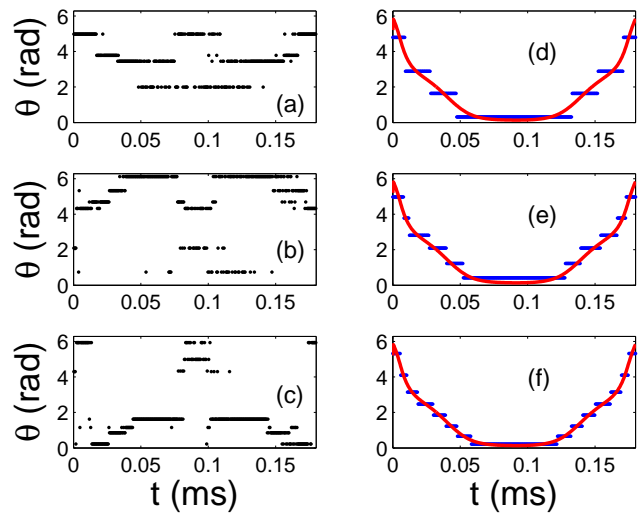}
\caption{(Color online) Examples of optimal discrete control phases obtained with the discrete GRAPE method (left column) and with the Lloyd's algorithm (right column)
for the cases $M=4$ (a-d), $M=6$ (b-e) and $M=8$ (c-f). For GRAPE, the initial discrete values are randomly chosen in the interval $[0, 2\pi[$. The continuous optimal field is plotted in solid line (red) in the right panels.}\label{fig5}
\end{figure}

As could be expected, the larger the number of discrete values $M$ is, the higher the efficiency. This point is clearly shown in Fig.~\ref{fig4} where the value of the final quality factor $\Phi$ is plotted for the two specific choices of the initial set of discrete values as a function of the number of discrete values $M$.

For the case of a random choice of the initial set, we observe that the averaging over 100 realizations of the optimized quality factor seems less efficient than the performance achieved when the initial discrete values are uniformly distributed in the interval $[0,2\pi[$. However, the maximum of the hundred realizations have a very high performance, which is better than the one with a uniform distribution, even for $M=4$. This point is illustrated in Fig.~\ref{fig4}, where the highest quality factor among the 100 random realizations is plotted as a function of $M$. This observation illustrates a standard optimization problem with a gradient procedure, which only gives a local information and cannot avoid traps and local maxima in the control landscapes \cite{reviewQC}.

The complexity of the optimal solutions designed by the algorithm is illustrated in Fig.~\ref{fig5} for $M=4,~6$ and 8. For these particular cases, we get an efficiency of the order of 0.99 for $M=4$, 0.9925 for $M=6$ and 0.9931 for $M=8$.

\begin{figure}[htpb]
\centering
\includegraphics[scale=0.75]{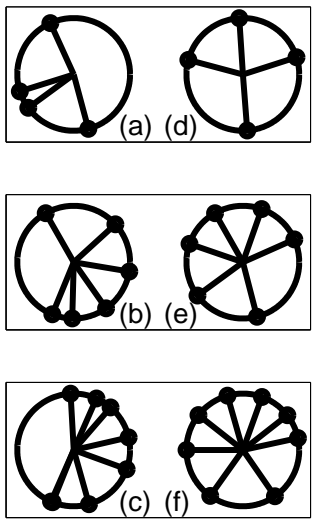}
\caption{Representation on a circle of the discrete phases computed in the example of Fig.~\ref{fig5} for the GRAPE (left column) and the Lloyd (right column) algorithms with $M=4$ (top, a-d), $M=6$ (middle, b-e) and $M=8$ (bottom, c-f).}\label{fig7}
\end{figure}

In a second series of simulations, we explore the efficiency of the Lloyd's algorithm for the same control problem. As explained in Sec. \ref{sec2c}, we first apply a standard continuous version of GRAPE. Figure~\ref{fig5} shows the different discrete fields designed by this procedure, which remain very close to the continuous one. Note the difference with respect to the solution determined from GRAPE in Fig.~\ref{fig5}. The representation used in Fig.~\ref{fig7} helps the visualization of the two different results.

\begin{figure}[htpb]
\includegraphics[scale=0.55]{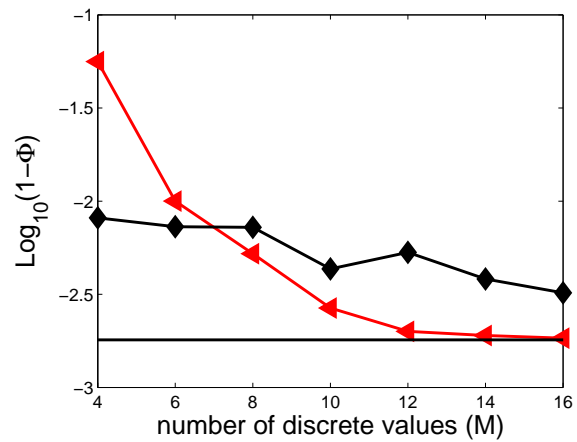}
\caption{(Color online) Comparison of the figures of merit achieved with the discrete version of GRAPE  and Lloyd's algorithm  for different numbers of discrete values. In the case of GRAPE, the best result over the 100 realizations of Fig.~\ref{fig2} has been used. The horizontal solid line indicates the efficiency of the continuous version of GRAPE.}\label{fig8}
\end{figure}

The numerical findings of Fig.~\ref{fig8} confirm the superiority of the discrete version of GRAPE over Lloyd's algorithm when $M\leq 6$. For larger values of $M$, we observe that the Lloyd's algorithm becomes more efficient. This point can be interpreted in the light of the local procedure used here to design the mapping. Optimizing such a mapping is intrinsically a global combinatorial optimization problem which cannot be easily replaced by a local search method. Another choice could be to combine a GRAPE approach, together with a heuristic global procedure such as genetic algorithms, which would be used to optimize the mapping \cite{reviewQC}. The global nature of the search method seems interesting even if the computational cost associated with this approach may be prohibitive.

As could be expected, the figure of merit for the Lloyd procedure tends to the continuous optimal result as the number of discrete values goes to infinity. This is not the case for the discrete version of GRAPE when random or uniform initial phases are used. In order to improve the efficiency of the quantization procedure, we have also used the values derived by Lloyd's algorithm as initial guess for GRAPE. However, we observe, for this example, that the combination of the two algorithms does not lead to a significantly better figure of merit. Again this is probably a direct consequence of the approach we choose to optimize the mapping. The comparison of the relative performance of the two algorithms should be tested on other control problems to confirm the conclusions of this paper.

\section{Conclusion}\label{sec4}
This paper has focused on the application of new optimization procedures for studying the simultaneous control of an ensemble of uncoupled spins with different offsets. A basic feature of the methods under consideration is that they allow us to account for a quantization of the control, which can only take a fixed finite number of values. Numerical tests have demonstrated the efficiency of the proposed procedures, showing that a high quality control can be achieved even with a very small number of quantized values. This conclusion is crucial in some electronically controlled nanodevices or in Nuclear Magnetic Resonance for which such constraints have to be taken into account. Experiments are in progress in the field of Electron Paramagnetic Resonance in which such constraints are imposed by the hardware \cite{EPR}.\\ \\
\noindent\textbf{ACKNOWLEDGMENT}\\
S.J. Glaser acknowledges support from the DFG (GI 203/7-2). M. Lapert acknowledges support from Alexander von Humboldt
Stiftung. D. Sugny and S. J. Glaser acknowledge support from the PICS program of
the CNRS and the ANR-DFG research program Explosys (ANR-14-CE35-0013-01; GL203/9-1). J.S was partially supported by the Agence Nationale
de la Recherche (ANR), Projet Blanc EMAQS number ANR-2011-BS01-017-01. This work has been done with the support of the Technische Universit\"at M\"unchen – Institute for Advanced
Study, funded by the German Excellence Initiative and the European Union Seventh Framework Programme under grant agreement 291763.

\bibliographystyle{apsrev}

\end{document}